\documentclass{amsart}[11pt]
\usepackage{myAMSpackages}
\usepackage{Universal}
\usepackage{MyAMSEnvironments}
\usepackage{graphicx,subfigure,epsfig}

\newcommand{\tor}{\text{Tor}}
\newcommand{\Tor}{\text{TOR}_2}

\newcommand{\tr}[1][H]{\ensuremath{\text{tr}_{#1}}\xspace}

\newcommand{\Ric}[1][H]{\ensuremath{\text{Rc}^{#1}}\xspace}
\newcommand{\srm}{\ensuremath{\text{Rm}^s}\xspace}

\newcommand{\rgv}{\widehat{\mathfrak{R}}}

\newcommand{\src}[1][]{\ensuremath{ \text{Rc}^s_{#1}}\xspace}

\newcommand{\lp}[1][H]{\triangle_{#1}}

\newcommand{\bp}{\begin{pmatrix}}
\newcommand{\ep}{\end{pmatrix}}

\begin{document}
\title[First eigenvalue of the horizontal Laplacian]{Bounds for the first eigenvalue of the horizontal Laplacian in positively curved sub-Riemannian manifolds.}
\author{Robert K. Hladky}
\address{North Dakota State University Dept. \#2750, 
PO Box 6050, 
Fargo ND 58108-6050}
\email{robert.hladky@ndsu.edu}

\begin{abstract}
We establish lower bounds for the first non-zero eigenvalue for the natural geometric sub-elliptic Laplacian operator defined on sub-Riemannian manifolds of step 2 that satisfy a positive curvature condition. The methods are very general and can be applied even when the sub-Riemannian geometry has considerable torsion.
\end{abstract}

\maketitle

\section{Introduction}\label{S:ID}
The classical Lichnerowicz theorem for Riemannian geometry states that for a Riemannian manifold $M^d$ with Ricci curvature $\text{Ric}\geq \rho>0$ then the first non-zero eigenvalue for the Laplacian is sharply bounded as follows
\bgE{L} \lambda_1 \geq \frac{ (d-1) \rho}{d} .\enE
The classical Myers theorem states that same condition on the Ricci curvature implies that $M$ is compact. These results illustrate the deep connections between analysis, geometry and topology on Riemannian manifolds.

In \cite{Greenleaf}, Greenleaf extended these results to strictly pseudoconvex manifolds of dimension $2d+1$ with $d \geq 3$. Working with the sub-Laplacian associated to the Tanaka-Webster connection, he was able to show that if
\bgE{psiEst} \text{Ric}(X,X) + \frac{d}{2} \aip{\tor(T,X)}{JX}{} \geq \rho \left| X \right|^2 \enE
for all horizontal $X$ then 
\bgE{spc}
\lambda_1 \geq \frac{  \rho d}{d+1}.
\enE
This was later extended to the case $d\geq 2$ by Li and Luk \cite{LiLuk}. Considerable effort has gone into studying the case $d=1$, see for example \cite{ChangChiu}. Additional assumptions are needed to reproduce \rfE{psiEst}, but there are general estimates that make use of more complicated bounds on the torsion.

In \cite{Hladky4}, the author introduced a notion of connection adapted  to sub-Riemannian manifolds that simultaneously generalizes the Lev-Civita and the Tanaka-Webster connections. Under some mild geometric assumptions,  there is a natural sub-elliptic Laplacian associated to this connection which is formally self-adjoint and negative.  The purpose of this paper is to study bounds of the first eigenvalue of this sub-elliptic operator on compact manifolds satisfying a positive curvature constraint. Since the torsion of the sub-Riemannian connection is typically complicated, its effect on the eigenvalue bounds is discussed in detail.  

The main result is the following:

\begin{thms}
If $M$ is a sub-Riemannian manifold with horizontal dimension $d$ such that there are constants $\rho_1,\rho_2>0$ and $0\leq x<1$ with\[  \mathcal{R}^x(A,A) \geq \rho_1 \left| A_H \right|^2 +\rho_2  \left| A_V\right|^2  \]
and
\[ \rho_1 >m(\omega,\chi,\psi) :=\inf\limits_{s>0} \left( s\omega + \frac{\chi}{s} +  \frac{\psi}{s^2}  \right)
 \] 
 then the smallest positive eigenvalue of $\lp$ satisfies the bound
\[  \lambda_1 \geq \dfrac{ \rho_1 - m(\omega,\chi,\psi) }{(1-x) \frac{d-1}{d}+\omega}.\] 
\end{thms}
Here $\omega,\chi$ and $\psi$ are  constants (invariant under constant vertical rescaling) that measure the torsion of $M$ and $\mathcal{R}^x(A,A)$ is a variation on the Baudoin-Garofalo tensor introduced in \cite{BaudoinGarofalo}, which should be thought of as a sub-Riemannian analogue to Ricci curvature.

The dependence of this eigenvalue bound on the constrained variable $x$ essentially reduces the problem of optimal estimates to a 1 dimensional maximization problem. These can be solved using elementary techniques, but general formulas are overly complicated.

Under certain conditions on the torsion which will be described in detail in Section 4, this result can be simplified greatly. The category of almost strictly normal manifolds is introduced. This category is large enough to include most examples traditionally studied, including strictly pseudoconvex pseudohermitian manifolds.

\begin{thms}
If $M$ is compact, almost strictly normal and there are  constants $0\leq x< 1$ and  $\rho_1,\rho_2>0$ such that
\[\begin{split} \mathcal{R}^x(A,A)  &+ 2 \aip{\tr \Tor(A_H)}{A_H}{} \geq 2 \left| \tau_V^H(A_H) \right| \, \left| \tau_H^V(A_H) \right|  \\ & \qquad  + \rho_1 \left| A_H\right|^2 + \rho_2\left|A_V \right|^2 \end{split}\]
then
\[ \lambda_1 \geq \dfrac{ \rho_1}{ (1-x) \frac{d}{d-1} + \omega}.\]
\end{thms}
Here $\tau_H^V$ and $\tau_V^H$ are semi-norms that again measure torsion, but have significantly simpler descriptions than the constants $\chi$ and $\psi$.

 In \rfS{EX}, several examples are produced to demonstrate the developed techniques and illustrate use of the dependence on $x$.

\section{Computation}\label{S:CP}

\bgD{sRC}
An sub-Riemannian manifold with complement (sRC-manifold) is a manifold $M$ together with smooth constant rank bundles $HM, VM$ such that 
\[ TM = HM \oplus VM\] 
and a smooth inner product on $HM$. 

A metric extension for $M$ is a Riemannian metric $g$ that agrees with the given inner product on $HM$ and $g(HM,VM)=0$.
\enD

The following theorem was shown in \cite{Hladky4}

\bgT{Connection}
If $M$ is an sRC-manifold with metric extension $g$ there is a unique connection $\nabla$ such that
\begin{itemize}
\item $\nabla g = 0$,
\item $\tor(HM,HM) \subseteq VM$, $\tor(VM,VM)\subseteq HM$,
\item For $X,Y \in HM$, $T,U \in VM$,
\begin{align*}
 \aip{ \tor(X,T)}{Y}{} &= \aip{X}{\tor(Y,T)}{}\\ \aip{ \tor(X,T)}{U}{} &= \aip{T}{\tor(X,U)}{}.\end{align*}
\end{itemize}
Furthermore if $X,Y$ are horizontal vector fields and $T$ is a vertical vector field then $\nabla X$, $\tor(X,Y)$, $\tor(X,T)_H$ are all independent of the choice of $g$.
\enT

These connections are not torsion-free and this presence of torsion greatly complicates analysis on sRC-manifolds as compared to the Riemannian case. To obtain and optimize results, we shall use a variety of restrictions on the torsion.

%\bgD{Constants}
%Let 
%\begin{align*}
 %\ub &= \frac{d}{2} \sup \left\{   \left|  \tor(T,X_1)_H \right|^2  \right\} \\
  %\ub_V &= \frac{d}{2} \sup \left\{   \left|  \tor(T,X_1)_V \right|^2  \right\} \\
 %\kappa &= \sup  \sum\limits_E \left|  \tor(E,X_1) \right|^2 \\
 %\sigma &=  \frac{1}{2} \sup   \left|\sum\limits_E \nabla  \aip{\tor(T,E,X)}{E}{} \right| \\
% \sigma &= \sup \tau(T,X) \\
 %\chi &= \sup \aip{T}{\tr (\nabla \tor)(T)}{}
 %\end{align*}
%\enD

\bgD{Normality}
Let $\{E_i\}$ be any local orthonormal frame for $HM$ and $\{U_\ua\}$ any local orthonormal frame for $VM$,
\begin{itemize}
\item An sRC-manifold is $H$-normal if $\tor(HM,VM) \subseteq VM$. This is independent of $g$.
\item A metric extension is $V$-normal if $\tor(HM,VM) \subseteq HM$.
\item A metric extension is strictly normal if $\tor(HM,VM)=0$.
\item  The rigidity tensor for $g$ is 
\[ \mathfrak{R}(A) = \sum\limits_k \aip{\tor(E_k,A)}{E_k}{} + \sum\limits_i \aip{\tor(U_i,A)}{U_i}{}. \]
 The rigidity vector for $g$ is 
 \[ \rgv = \sum\limits_k \mathfrak{R}(E_k)E_k + \sum\limits_i \mathfrak{R}(U_i)U_i .\]
 The sRC-manifold is $H$-rigid if $\rgv_H \equiv 0$. The metric extension is $V$-rigid if $\rgv_V \equiv 0$ and totally rigid if $\rgv \equiv 0$.
\end{itemize}
\enD

To measure the size of the torsion on $M$, we introduce the following metric extension dependent semi-norms: 
\begin{align}
\left|\tau_V^H(A) \right|^2  &=  \sum\limits_{i,k} \aip{\tor(A,U_k)}{E_i}{}^2,  \\
\left| \tau_H^V(A)\right|^2 & = \sum\limits_{i,k} \aip{\tor(A,E_i)}{U_k}{}^2, \\
\left| \tau_H(A) \right|^2 &= \sum\limits_{i,j} \aip{\tor(E_i,E_j)}{A}{}^2
\end{align}
where $E_i$ is an orthonormal frame for $HM$ and $U_k$ is an orthonormal frame for $VM$.

We shall also need a sub-Riemmanian equivalent of the Ricci curvature. At first glance, the most natural generalization to sRC-manifolds would appear to be
\[ \tr \srm(A,B) = \sum\limits_k \srm(E_k,A,B,E_k), \]
where $E_k$ is any orthonormal frame for $HM$. However, this tensor is not in general symmetric  when restricted to horizontal vectors. Part of the reason for this is due to second order torsion terms, which of course are absent in the Riemannian setting. We introduce the tensor
\bgE{tor2}
\Tor(A,B,C) = \tor(A,\tor(B,C))
\enE
and using this we can define the following.

\bgD{sRicci} The sub-Ricci curvature of an sRC-manifold $M$ is the tensor
\begin{align*}
 \src[](A,B) &= \tr \srm(A,B) -  \frac{1}{2}  \sum\limits_k  \aip{\Tor(E_k,A_H,B_H)}{E_k}{}\\
 & \qquad -  \aip{\tr \Tor(A_H)}{B_H}{}\\
\end{align*}
where $\{E_k\}$ is any horizontal orthonormal frame. This tensor is independent of the choice of metric extension. 
\enD

\bgL{asym}
For $X,Y \in HM$ and $T \in VM$, the sub-Ricci curvature satisfies
\begin{align*}
\src (X,Y) &= \src (Y,X), \\
\src( X, T) & =0, \\
\src (T, X) &= \sum\limits_k \aip{ (\nabla \tor  - \Tor )(E_k,T,X) }{E_k}{}  - \aip{\Tor(T,X,E_k)}{E_k}{} \\
& \qquad + \aip{\tr (\nabla \tor - \Tor) (T)}{X}{}.
 \end{align*}
  \enL

\pf The middle result is trivial. To show the others, we apply a  standard result in differential geometry (see \cite{Hladky4} for detail), that for any curvature 
\[ \mathcal{C}R(A,B)C = -\mathcal{C}\Tor(A,B,C) + \mathcal{C} (\nabla \tor)(A,B,C) \]
where $\mathcal{C}$ represents cyclic permutation. Then following for example Proposition 7.4 in \cite{Lee3}, we see that from elementary properties of  curvature that
\begin{align*}
  2 \aip{R(C,A)B}{D}{} &- 2 \aip{R(B,D)C}{A}{} = \mathcal{C} \aip{\mathcal{C}R(A,B)C}{D}{}  \\
&= -\mathcal{C}\aip{ \mathcal{C}\Tor(A,B,C)}{D}{} + \mathcal{C}\aip{\mathcal{C} (\nabla \tor)(A,B,C)}{D}{}.
\end{align*}
First we apply this with $C=D=E \in HM$ and $A=X$, $B=Y$, to see that 
\begin{align*}  \aip{ R(E,X)Y}{E}{} & - \aip{R(E,Y)X}{E}{} = \aip{ \Tor(E,X,Y)}{E}{}\\
& \qquad + \aip{\Tor(E,E,X)}{Y}{}  - \aip{\Tor(E,E,Y)}{X}{}. 
\end{align*}
Secondly, since $ \aip{R(E,X)T}{E}{}=0$, if instead we set $C=D=E$, $A=T$, $B=X$, we have
\begin{align*}
2 \aip{R(E,T)X}{E}{}  &=  \mathcal{C}\aip{ \mathcal{C}(\nabla \tor)(T,X,E)}{E}{}  - \mathcal{C} \aip {\mathcal{C} \Tor(T,X,E) }{E}{}\\
&= 2 \aip{ \nabla \tor (E,T,X)}{E}{} +2 \aip{\nabla \tor (T,E,E)}{X}{} \\ 
& \qquad -2 \aip{ \Tor(E,E,T) }{X}{} -2 \aip{\Tor(E,T,X)}{E}{} \\
& \qquad  -2\aip{\Tor(T,X,E)}{E}{}.
\end{align*}
Next we sum over a frame for $HM$ to yield the remaining results.
 
\epf

\bgD{HL}
For a tensor $\tau$, the horizontal gradient of $\tau$ is defined by
\[ \nabla_H \tau = \nabla_{E_i} \tau \otimes E_i, \]
where $E_i$ is a (local) orthonormal frame for $HM$. 
The horizontal Hessian of $\tau$ is defined by
\[ \nabla^2_H \tau (X,Y)=\left(  \nabla_{X} \nabla_{Y} - \nabla_{\nabla_{X}Y} \right)  \tau\]
for $X,Y \in HM$ and zero if either $X,Y \in VM$.
Finally, the horizontal Laplacian of $\tau$ is defined by
\[ \lp \tau = \tr \left( \nabla^2_H \tau \right)  =\sum_i \left(  \nabla_{E_i} \nabla_{E_i} - \nabla_{\nabla_{E_i}E_i} \right)  \tau.  \]
\enD

The Laplacian on a Riemannian manifold has a rich and interesting $L^2$-theory. To replicate this for sRC-manifolds, it is necessary to choose a metric extension. This metric extension then yields a volume form and we have meaningful $L^2$-adjoints. Unfortunately, the horizontal Laplacian defined here, does not always behave as nicely as the Riemannian operator. However, if we make a mild assumption on the metric extension, much of the theory can be generalized.  The following was shown in \cite{Hladky4}

\bgL{VR}
Suppose that $M$ is orientable and has a $V$-rigid metric extension $g$.  Then on functions,
\[ \lp  = E_i^2 + \text{div } E_i =- \nabla^*_H \nabla_H  \]
where the divergence and $L^2$ adjoint are taken with respect to the volume form $dV_g$.
\enL
Thus on a $V$-rigid sRC-manifold, the horizontal Laplacian behaves qualitatively in a similar fashion to the Riemannian Laplacian. If $M$ does not admit a $V$-rigid extension, then the horizontal Laplacian becomes substantially harder to work with. In fact, we shall often need to assume the stronger condition that the metric extension is $V$-normal.

A key result from \cite{Hladky4} is the following Bochner theorem

\bgT{Bochner} If $F=\nabla f$ and $\bullet \in \{H,V\}$
\[ \begin{split} \frac{1}{2} \lp  \left| F_\bullet \right|^2 &- \aip{\nabla_\bullet \lp f}{F_\bullet}{} =\tr \srm (F_\bullet,F_H)  + \left\| \nabla^2_{H} F_\bullet \right\|^2\\ 
& \qquad -  2S(F,F_\bullet)  + \aip{F}{ \tr (\nabla \tor) (F_\bullet) }{}  \\
& \qquad - \aip{F}{\tr \Tor(F_\bullet)}{} \\
\end{split}\]
where $E_i$ is any local orthonormal frame for $HM$ and
\begin{align*}
\tr (\nabla \tor)(A) &= \sum\limits_i \nabla \tor(A,E_i,E_i),   \\
\tr \Tor(A) &= \sum\limits_i \Tor(E_i,E_i,A), \\
S(A,B) &= \sum\limits_i  \aip{\nabla_{E_i} A}{\tor(E_i,B)}{}.
\end{align*}
\enT

We now introduce the following tensors and forms
\bgD{Forms}
For $0\leq x<1$, the BG-curvatures $\mathcal{R}^x$ are the unique  2-tensors such that
\[ \begin{split} \mathcal{R}^x(A,A) &= (1-x) \src (A_H,A_H) +(1+x) \aip{\tr (\nabla \tor) A_H}{A}{}   \\ & \qquad + \frac{1+3x}{4} \left| \tau_H(A) \right|^2 -(1-x) \aip{A_V}{\tr \Tor(A_H)}{}.\end{split}\]
 We remark briefly that, if $M$ is strictly normal, $\mathcal{R}^0$ agrees with the Baudoin-Garofalo curvature tensor introduced in \cite{BaudoinGarofalo}.

Following \cite{BaudoinGarofalo}, we define symmetric bilinear forms  for $\bullet \in \{H,V\}$ by  
\begin{align*}
\Gamma^\bullet(f,g) &= \aip{\nabla_\bullet f }{\nabla_\bullet g}{} ,\\
\Gamma^\bullet_2(f,g) &= \frac{1}{2} \lp \Gamma^\bullet(f,g)  - \Gamma^\bullet(\lp f,g) - \Gamma^\bullet(f, \lp g). \\
\end{align*}
\enD

For any symmetric bilinear form, we shall adopt the notation $B(f) = B(f,f)$, so for example
\[\Gamma^H(f)+\Gamma^V(f) = \| \nabla f\|^2 .\]

It was also shown in \cite{Hladky4} that
\[ \|\nabla^2_H f\|^2 = \|\nabla^{2,sym}_H f\|^2+ \frac{1}{4} \left| \tau_H(F) \right|^2.\]
We shall further decompose the symmetric part as
\bgE{spl}
\begin{split}
 \| \nabla^{2,sym}_{D} f \|^2&= \sum\limits_i \left| \nabla^{2,sym} f(E_i,E_i) \right|^2,\\
   \|\nabla^{2,sym}_{*}  f\|^2&= \sum\limits_{i\ne j} \left| \nabla^{2,sym} f(E_i,E_j) \right|^2
   \end{split}
   \enE
and note that, in particular
\bgE{HessLp} \| \nabla^{2,sym}_{D} f \|^2 \geq \frac{1}{d} \left| \sum\limits_i  \nabla^{2,sym} f(E_i,E_i) \right|^2 = \frac{1}{d} \left( \lp f \right)^2.\enE

The purpose of  seemingly complicated tensors is the following integrated horizontal Bochner formulas
\bgT{HBochner}
For $0\leq x \leq 1$ and $\nu>0$
 \begin{align*}
 \int (1-x) \Gamma^H_2(f)   & \geq  \int  \mathcal{R}^x(F_H,F_H) + \frac{1-x}{d} \left( \lp f \right)^2 \; dV\\
 &  \qquad  - \int  \frac{1}{\nu} \left| \tau_H^V(F_H) \right|^2  +  \nu \left| \nabla_H F_V \right|^2  \; dV. \end{align*}
   \enT

\pf  The pointwise Bochner identity of \rfT{Bochner} can be written 
\bgE{GCI}
\begin{split}
\Gamma^H_2(f) &= \mathcal{R}^0(F,F) + \|\nabla^{2,sym}_H f\|^2 -2S(F,F_H).
\end{split}
\enE 

Now we compute that 
\begin{align*}
-2&\int S(F,F_H) dV  = -2  \sum\limits_k \int   \aip{\nabla_{E_k} F_V}{\tor(E_j,F_H)}{} dV \\
&= 2  \sum\limits_k \int \aip{F}{\nabla \tor(E_k,F_H,E_k) + \tor(E_k,\nabla_{E_k} F_H)}{} dV \\
&= \int -2 \aip{F}{\tr (\nabla \tor)(F_H)}{} + 2 \sum\limits_{i,j} \aip{F}{\tor(E_i,E_j)}{}\aip{\nabla_{E_i} F_H}{E_j}{}  dV\\
&=  \int  \sum\limits_{i,j}   \aip{F}{\tor(E_i,E_j)}{}\left( \aip{\nabla_{E_i} F_H}{E_j}{}+\aip{\nabla_{E_j} F_H}{E_i}{}   \right)  \\ & \qquad + \sum\limits_{i,j}   \aip{F}{\tor(E_i,E_j)}{}\aip{F}{\tor(E_j,E_i)}{}   \\
& \qquad  -2 \aip{F}{\tr (\nabla \tor)(F_H)}{}  dV \\
&= - \int  \left| \tau_H(F) \right|^2-2 \aip{F}{\tr (\nabla \tor)(F_H)}{}  dV
\end{align*}
as the first term of the penultimate line vanishes due to a symmetry/skew-symmetry argument. It should be remarked here, that as the frame $E_i$ is only defined locally, that we should technically employ a partition of unity argument to conduct this integration-by-parts. All terms are independent of the choice of horizontal frame, so this argument is routine.

From this we see that \rfE{GCI} can be integrated to 
\bgE{BochN}
\begin{split}
 \int \Gamma^H_2(f) dV &= \int \mathcal{R}^0(F,F) + \| \nabla_{D}^{2,sym} f\|^2 +\|\nabla^{2,sym}_*f\|^2  \\
& \qquad  -2 \aip{F}{\tr (\nabla \tor)(F_H)}{}- \left| \tau_H(F)\right|^2  dV.
\end{split}
\enE
In particular,
\bgE{offest}
\begin{split}
 \int \| \nabla_{*}^{2,sym} f\|^2 &= \int \Gamma^H_2(f) - \mathcal{R}^0(F,F)  - \| \nabla^{2,sym}_D f(E_i,E_i)\|^2 \\
& \qquad  + 2 \aip{F}{\tr (\nabla \tor)(F_H)}{} + \left| \tau_H(F)\right|^2  dV
\end{split}
\enE
So for $0\leq x \leq 1$
\begin{align*}
\int \Gamma^H_2 (f) &  dV  \geq  \int x \Gamma^H_2(f) + (1-x) \mathcal{R}^0(F,F) + (1-x)\|\nabla_D^{2,sym} f\|^2 \\
& \qquad  +2x  \aip{F}{\tr (\nabla \tor)(F_H)}{}  + x \left| \tau_H(F)\right|^2 - 2S(F,F_H) dV \\
& \geq  \int x \Gamma^H(f) + \mathcal{R}^x(F,F)  + (1-x) \left( \lp f \right)^2 - 2S(F,F_H)  dV.
\end{align*} 
The result then follows from the observation that for all $\nu >0$
\[ 2 \left| S(F,F_H) \right| \leq  \nu \left| \nabla_H F_V \right|^2 + \frac{1}{\nu} \left| \tau_H^V(F_H) \right|^2.\] \epf

The problematic term in the last equation is $\left| \nabla_H F_V \right|^2$, which thus far cannot be controlled. Accordingly, we now turn our attention $\Gamma^V_2$ and the vertical Bochner formula. The key result is the next theorem. 
  
 \bgT{VBochner} If $F = \nabla f$ then
\[  \int_M \Gamma^V_2(f) dV  = \int_M  \| \nabla_H F_V \|^2 - \mathcal{T}_1(F_H,F_V) - \mathcal{T}_2(F_H,F_H)  \; dV \]
 where  mixed distortion tensor $T_1$ is defined by
\begin{align*}
 \mathcal{T}_1(A,B) &=  \sum\limits_k \Big( \aip{ (\Tor - \nabla \tor )(E_k,A_V,B_H) }{E_k}{}   \\ & \qquad +\aip{\Tor(A_V,B_H,E_k)}{E_k}{} \Big) \\
& \qquad + 4 \aip{\tr \Tor(A_V)}{B_H}{}
  \end{align*}
and the pure distortion tensor $T_2$ is the unique symmetric 2-tensor with
  \begin{align*}
   \mathcal{T}_2(A,A) &=     2\left| \tau_V^H(A_H) \right|^2  +  \aip{ \tr[V](\nabla \tor)(A_H)}{A_H}{} \\
&  \qquad  + \aip{\tor(A_H,\rgv_V)}{A_H}{}. \\
  \end{align*}
\enT

The proof of this theorem follows trivially from \rfL{asym}, \rfT{Bochner} and the following lemma: 

\bgL{et} 
The error term $S(F,F_V)$ integrated as either
\[ \begin{split} 2 \int_M&  S(F,F_V)  dV =  \int_M \aip{F_V}{ \tr \nabla \tor(F_V)}{}  \\ 
& \qquad +2\left| \tau_V^H(F_H) \right|^2   +2\aip{F_H}{\tr (\nabla \tor+\Tor)(F_V)}{}   \\
&  \qquad \qquad -2 \sum\limits_{i}  \aip{\tor(F_H,U_i)}{\nabla_{U_i} F_H}{} \; dV\\
 &=  \int_M \Big(  \aip{F_V}{ \tr (\nabla \tor)(F_V)}{}  +2\aip{F_H}{\tr (\nabla \tor + \Tor)(F_V) }{}   \\
&  \qquad + 2\left| \tau_V^H(F_H) \right|^2   \\
&  \qquad +  \aip{ \tr[V](\nabla \tor)(F_H)}{F_H}{}  + \aip{\tor(F_H,\rgv_V)}{F_H}{}  \Big]\Big) \; dV 
\end{split}.\]
\enL

\pf

First we note
\bgE{bk}
 \aip{\nabla_E F}{\tor(E,F_V)}{} =  \aip{\nabla_E F_H}{\tor(E,F_V)}{} +  \aip{\nabla_E F_V}{\tor(E,F_V)}{}.
\enE
Integrating by parts (again suppressing a partition of unity decomposition) yields 
 \begin{align*}
 \int_M  & \sum\limits_k   \aip{\nabla_{E_k} F_V}{\tor(E_k,F_V)}{} dV \\
 & = \int_M  \sum\limits_{k} \Big(  {E_k} \aip{F_V}{\tor(E_k,F_V)}{} -  \aip{F_V}{\nabla_{E_K} \tor(E_k,F_V)}{} \Big) dV \\
& =  -  \int_M  \sum\limits_k \Big( \aip{F_V}{\nabla \tor(E_k,F_V,E_k) }{} + \aip{F_V}{ \tor(E_k,\nabla_{E_k} F_V)}{}  \\ 
&  \qquad  + \aip{F_V}{\tor( \nabla_{E_k} E_k,F_V)}{}  +  \text{div} E_k   \aip{F_V}{\tor(E_k,F_V)}{}   \Big)  dV\\
 &= - \int_M \sum \limits_k \Big( \aip{F_V}{\nabla \tor(E_k,F_V,E_k) }{} + \aip{F_V}{ \tor(E_k,\nabla_{E_k} F_V)}{} \Big) dV.\end{align*}
However we can now use torsion symmetry to see that
\bgE{et1} 2 \int_M \sum\limits_k \aip{\nabla_{E_k} F_V}{\tor(E_k,F_V)}{} dV =  \int_M \aip{F_V}{ \tr \nabla \tor(F_V)}{} dV. \enE

To deal with the other term of  \rfE{bk} we apply a similar integration-by-parts argument to see 
\bgE{et2}
\begin{split}
 \int_M & \sum\limits_k   \aip{\nabla_{E_k} F_H}{\tor(E_k,F_V)}{} dV  \\
&   = -  \int_M \sum\limits_k \Big(  \aip{F_H}{\nabla \tor(E_k,F_V,E_k)  }{}   + \aip{F_H}{ \tor(E_k,\nabla_{E_k} F_V) }{}  \Big)dV\\
 &  =  \int_M \aip{F_H}{\tr (\nabla \tor)(F_V) }{}dV- \int_M \sum\limits_k \aip{F_H}{\tor(E_k,\nabla_{E_k} F_V)}{}.
\end{split}\enE

Now since $F = \nabla  f$ is a closed vector field we see that
\bgE{closed} \aip{\nabla_A F}{B}{} - \aip{\nabla_B F}{A}{}= \aip{\tor(B,A)}{F}{}\enE
and so
\begin{align*}
\sum\limits_k & \aip{ \tor(E_k,\nabla_{E_k} F_V)}{F_H}{} =\sum\limits_{i,k}  \aip{\tor(E_k,U_i)}{F_H}{}\aip{\nabla_{E_k} F_V}{U_i}{} \\ 
&  =\sum\limits_{i,k} \aip{\tor(E_k,U_i)}{F_H}{} \left( \aip{\nabla_{U_i} F_H}{E_k}{} + \aip{F}{\tor(U_i,E_k)}{} \right)\\
& =\sum\limits_{i,k} \Big(  - \aip{\tor(E_k,U_i)}{F_H}{}^2 -\aip{\tor(E_k,U_i)}{F_H}{}\aip{\tor(E_k,U_i)}{F_V}{}  \\
&  \qquad +  \aip{\tor(E_k,U_i)}{F_H}{} \aip{\nabla_{U_i} F_H}{E_k}{} \Big)  \\
& = -\left| \tau_V^H(F_H) \right|^2   -\aip{F_H}{\tr \Tor(F_V)}{}  + \sum\limits_{i}  \aip{\tor(F_H,U_i)}{\nabla_{U_i} F_H}{}.
\end{align*}
Integrating-by-parts the yields
\begin{align*}
\int_M  U_i \aip{\tor(X,U_i)}{X}{}  &- \aip{\tor(X,\nabla_{U_i} U_i)}{X}{} dV  \\
 &=\int_M -\text{div}_g U_i \aip{\tor(X,U_i)}{X}{}  \\
 & \qquad + \aip{\nabla_{U_j} U_i}{U_j}{}  \aip{ \tor(X,U_i)}{X}{} dV \\
&=\int_M - \mathfrak{R}(U_i )\aip{\tor(X,U_i)}{X}{} dV.
\end{align*}
 But computing directly
 \begin{align*}
  U_i \aip{\tor(X,U_i)}{X}{}  &- \aip{\tor(X,\nabla_{U_i} U_i)}{X}{} \\
  &= \aip{ \nabla_{U_i} \tor(X,U_i)}{X}{} + \aip{\tor(X,U_i)}{\nabla_{U_i} X}{}  \\ & \qquad - \aip{\tor(X,\nabla_{U_i} U_i)}{X}{} \\
  &= 2 \aip{\tor(X, U_i)}{\nabla_{U_i} X}{} + \aip{ \nabla \tor(X,U_i,U_i)}{X}{}.
 \end{align*}
 Comparing we conclude
 \bgE{et3}
 \begin{split}
 2 \int_M &   \aip{\tor(F_H , U_i)}{\nabla_{U_i} F_H}{}  \\ & =  - \int_M \sum\limits_{i}  \aip{ \nabla \tor(F_H,U_i,U_i)}{X}{}  +  \mathfrak{R}(U_i )\aip{\tor(F_H,U_i)}{X}{} dV.
 \end{split}
\enE

%Combining \rfE{et1}, \rfE{et2} and \rfE{et3} yields the result
%\begin{align*}
%- 2\int_M \aip{\tor(E,\nabla_E F_V)}{F_H}{} &= 2 \int_M \left| \tor(F_H,U)\right|^2  + \aip{ \tr[V] (\nabla \tor)(F_H)}{F_H}{}  \\
%& \qquad + \aip{\tr \Tor(F_V)}{F_H}{}  \\
%& \qquad + \sum\limits_{i,k} \aip{E_k}{\tor(E_k,U_i)}{} \aip{\tor(F_H,U_i)}{F_H}{} dV
%\end{align*}
Combining all of these yields the result.

\epf

From \rfT{HBochner} and \rfT{VBochner} we obtain the following integrated curvature-dimension inequality: for any $\nu >0$, 

\bgE{icdi}
\begin{split}
\int   (1-x)\Gamma^H_2(f)  &+ \nu \Gamma^V_2(f) dV \\
&  \geq  \int \Big(  \mathcal{R}^x(F,F)  + \frac{1-x}{d} \left( \lp f\right)^2  - \frac{1}{\nu} \left| \tau_H^V(F_H) \right|^2   \\
&  \qquad  - \nu \mathcal{T}_1(F_H,F_V) - \nu \mathcal{T}_2 (F_H) \Big)\;  dV.
\end{split}
\enE 

\section{General eigenvalue estimates}

With the integrated curavture-dimension inequality \rfE{icdi} under our belts, we can now turn our attention to obtaining general eigenvalue bounds.  Our fundamental assumption will be that the sRC-manifold with metric extension is positively curved in the sense that there exist constants $0\leq x <1$ and $\rho_1,\rho_2 >0$ such that

 \bgE{A}
 \begin{split}
  \mathcal{R}^x(A,A) &\geq \rho_1 \left| A_H \right|^2 +\rho_2  \left| A_V\right|^2  \\
\end{split}
\enE 

The constant $\rho_1$ is invariant under constant vertical conformal rescaling, but $\rho_2$ is not. To address this, we also introduce
\bgE{omega} \kappa = \sup\limits_{X \in HM \backslash\{0\} }  \sum\limits_{k}  \dfrac{ \left| \tau_H^V(X)  \right|^2 }{|X|^2}, \qquad \omega = \frac{\kappa}{\rho_2}  \enE
Note that $\omega$ is an invariant under constant vertical conformal rescaling. To control the pure distortion term we  define a scale invariant constant $\chi \in \rn{}$ by 
\[ \chi =   \rho_2  \sup \limits_{X \in HM \backslash\{0\} } \dfrac{\mathcal{T}_2(X,X)}{\|X\|^2}\] 
In particular \[  \mathcal{T}_2(F_H,F_H) \leq \frac{\chi}{\rho_2}  \Gamma^H(f) .\]
To control the mixed distortion term, we introduce the constants
\[ \sigma = \sup\limits_{X \in HM \backslash\{0\}, T \in VM \backslash\{0\} }  \dfrac{ \mathcal{T}_1(T,X) }{|T| \, |X|}, \qquad \psi = \rho_2 \sigma^2. \]
 Again $\sigma$ is not scale invariant, but $\psi$ is. In particular, we have that for all $s>0$
\[ \mathcal{T}_1(F_V,F_H) \leq 2 \sigma \left| F_V \right| \; \left| F_H \right| \leq  \frac{\sigma^2}{s} \Gamma^H(f)  +  s\Gamma^V(f)\]

Now from \rfE{icdi} and assumption \rfE{A} we get that for any $\nu>0$
\bgE{icdi2}
\begin{split}
\int_M (1-x)\Gamma_2^H(f)  & +\nu  \Gamma_2^V(f) \; dV \geq  \int_M \Big[ \frac{1-x}{d} (\lp f)^2+   \left( \rho_1 - \frac{\kappa}{\nu} \right) \Gamma^H(f)  \\
& \qquad   + \rho_2 \Gamma^V(f)  - \nu \mathcal{T}_1(F_H,F_V) - \nu \mathcal{T}_2(F_H,F_H)  \Big] \; dV
\end{split}
\enE

Now if we choose $f$ such that $\lp f = -\lambda f$ with $\lambda >0$ and rescale so that $\int_M \Gamma^H(f) dV = 1$, then a simple integration by parts of the left hand side of \rfE{icdi} produces
\bgE{lcdi}
\begin{split}
(1-x) \frac{d-1}{d} \lambda^2  &\geq \left( \rho_1 - \frac{\kappa}{\nu}\right) \lambda + \int_M  \left( \rho_2 -\lambda \nu\right) \Gamma^V(f) dV \\
& \qquad  \qquad - \nu \int_M  \mathcal{T}_1(F_V,F_H) +  \mathcal{T}_2(F_H,F_H)  \; dV \\
&  \geq \left( \rho_1 - \frac{\kappa}{\nu} - \frac{\nu \sigma^2}{s} - \frac{\chi \nu}{\rho_2}  \right) \lambda  + \int_M  \left( \rho_2 -( \lambda +s) \nu\right) \Gamma^V(f)  \; dV
\end{split}
\enE
Setting $\nu = \frac{\rho_2}{\lambda + s}$ then yields
\bgE{B}
\begin{split}
 \left((1-x) \frac{d-1}{d} +\omega\right) \lambda  &\geq  \rho_1 - s \omega- \frac{\psi}{s(\lambda+s)} -\frac{\chi}{\lambda+s} \end{split}
 \enE
Thus if we set 
\bgE{m}
m(\omega,\chi,\psi) = \inf\limits_{s>0} \left( s\omega + \frac{\chi}{s} +  \frac{\psi}{s^2}  \right)
\enE
then \rfE{B} implies
\bgE{C}
\left((1-x)\frac{d-1}{d}+\omega\right) \lambda \geq \rho_1 - m(\omega,\chi,\psi)
\enE

\bgT{MainEst}
If $M$ satisfies the fundamental assumption \rfE{A} with
\[ \rho_1 > m(\omega,\chi,\psi)\]
then the smallest positive eigenvalue of $\lp$ satisfies the bound
\[  \lambda_1 \geq \dfrac{ \rho_1 - m(\omega,\chi,\psi) }{(1-x) \frac{d-1}{d}+\omega}.\]
Furthermore if
\begin{itemize}
\item If $\mathcal{T}_1 \equiv\mathcal{T}_2 \equiv 0 $  then \[m(\omega,\chi,\psi) = 0\]
\item  If $\mathcal{T}_1 \equiv 0$ then 
\[ m(\omega,\chi,\psi)  = 2\sqrt{ \omega\chi} \]
\item If $\mathcal{T}_2 \equiv 0 $  then \[m(\omega,\chi,\psi) =  \sqrt[3]{\frac{27 \omega^2\psi}{4}}\]
\end{itemize}
\enT

In the special case where $\mathcal{T}_1\equiv 0$, we can improve on this  result. To compress notation, set
\[ \Delta= (1-x) \frac{d-1}{d} .\]Since $\psi=0$, we can let $s\to 0$ in \rfE{B} to obtain
\[ \left(\Delta+ \omega\right) \lambda \geq  \rho_1 - \frac{ \chi }{\lambda} \]
or equivalently
\[ \left(\Delta+ \omega\right) \lambda^2 - \rho_1 \lambda + \chi  \geq 0 \]
However, we now note that we can apply this result with $\rho_2$ replaced $t \rho_2$ for any scalar $0<t\leq 1$.  This will rescale $\chi$ and $\omega$ to $\chi t$, $\frac{\omega}{t}$ respectively. Thus for all $0 < t \leq 1$ we have
\[  \left(\Delta+ \frac{\omega}{t}\right) \lambda^2 - \rho_1 \lambda +t \chi  \geq 0 \]
Thus if
\[ \rho_1^2 \geq 4 \omega \chi  + 4 \chi t  \frac{d-1}{d}  \]
Then
\[ \lambda  \notin \left( \;  \dfrac{ \rho_1 - \sqrt{\rho_1^2 - 4 \chi t \Delta - 4\omega \chi }}{2 (\Delta+\frac{\omega}{t})} \, , \,  \dfrac{ \rho_1 + \sqrt{\rho_1^2 -  4\chi t \Delta - 4\omega \chi  }}{2 (\Delta+\frac{\omega}{t})} \; \right). \]

Now as $t \to 0$, the left hand limit of the prohibited interval tends to $0$ also. Thus we have prohibited $\lambda$ from lying in a family of overlapping intervals, whose union is of the form $(0,b)$. Hence we have established the following theorem

\bgT{T1z}
If $M$ satisfies  $\mathcal{T}_1 \equiv 0$ and the fundamental assumption  \rfE{A} with
\[ \rho_1^2 > 4 \omega \chi \]
then for all $0<t\leq  \min \left\{ 1, \frac{\rho_1^2 - 4 \omega \chi}{4\Delta \chi} \right\}$
\[ \lambda_1 \geq  \dfrac{ \rho_1 + \sqrt{\rho_1^2 -  4 \Delta \chi t - 4\omega \chi  }}{2 (\frac{d-1}{d}+\frac{\omega}{t})}\]
In particular, 
\begin{enumerate}
\item If $0<\rho_1^2 - 4\omega \chi < 4 \chi \Delta$ then
\[ \lambda_1 \geq \dfrac{ \rho_1}{  \Delta \left( 2 +  \frac{4\chi}{\rho_1^2-4\omega \chi} \right) } \]
\item If $\rho_1^2 - 4\omega \chi \geq  4 \chi \Delta$ then
\[ \lambda_1 \geq  \dfrac{ \rho_1 + \sqrt{\rho_1^2 -  4 \Delta \chi  - 4\omega \chi  }}{2 (\Delta+\omega)}\]
\end{enumerate}
\enT

We conclude this section with the following easily verifiable remarks concerning the distortion tensors. 
\bgL{dist} \hfill

\begin{enumerate}
\item If $M$ is $H$-normal then $\mathcal{T}_2\equiv 0$.
\item If $M$ is $H$-normal and $VM$ is integrable then $\mathcal{T}_1 \equiv 0$.
\item If $M$ is strictly normal and $VM$ is integrable then $\mathcal{T}_1 \equiv 0 \equiv \mathcal{T}_2$.
\end{enumerate}
\enL

\section{Almost strictly normal manifolds}

\bgD{asn}
An sRC-manifold $M$ with metric extension $g$ is almost strictly normal if $M$ is $H$-rigid, $VM$ is integrable and the extension is $V$-normal.
\enD

\bgL{HrN}
If $M$ is $H$-rigid then
\[ \sum\limits_i \aip{\nabla \tor(E_i,T,X)}{E_i}{}=0\]  for all $X \in HM$ and $T \in VM$.
\enL

\pf 
Suppose $M$ is $H$-rigid, then \[ \sum \aip{\tor(E_i,T)}{E_i}{}=0\]
for all $T \in VM$. Differentiating produces
\begin{align*}
X \aip{\tor(E_i,T)}{E_i}{} &=  \aip{\tor(E_i,T)}{\nabla_X E_i}{} + \aip{\nabla_X \tor(E_i,T)}{E_i}{}\\
&= 2 \aip{\tor(E_i,T)}{\nabla_X E_i}{} + \aip{\nabla \tor(E_i,T,X)}{E_i}{}  \\ & \qquad + \aip{\tor(E_i,\nabla_X T}{E_i}{}\\
&= 2 \aip{\tor(E_i,T)}{E_j}{}\aip{E_j}{\nabla_X E_i}{} + \aip{\nabla \tor(E_i,T,X)}{E_i}{}  \\ & \qquad + \aip{\tor(E_i,\nabla_X T}{E_i}{}
\end{align*}
Since the first term vanishes by a symmetry/skew-symmetry argument and last term vanishes due to $H$-rigidity, the result follows immediately.

\epf

\bgC{HrN}
For an almost strictly normal manifold $\mathcal{T}_1 \equiv 0$ and
\[ \mathcal{T}_2(A,A) = 2 \left| \tau_V^H(A_H)\right|^2 + \aip{ \tr[V](\nabla \tor)(A_H)}{H_H}.\]
\enC

We shall also need the following: 
\bgL{VV}
If $T \in VM$ then 
\[ \aip{T}{\tr \Tor(T)}{} = \left| \tau_H^V(T) \right|^2. \]
\enL

\pf  This follows from direct computation.
\begin{align*}
\aip{T}{\tr \Tor(T)}{} &= \aip{T}{\tor(E_j,E_k)}{} \aip{\tor(E_j,T)}{E_k}{}  \\& \qquad + \aip{T}{\tor(E_j,U_i)}{}\aip{\tor(E_j,T)}{U_i}{} \\
&= \left| \tor(E_j,T) \right|^2 = \aip{\tor(E_j,U_i)}{T}{}^2.
\end{align*}
\epf

\bgL{tr}
If $(M,g)$ is almost strictly normal
\[ \int \Gamma_2^V(f)  \; dV= \int \left| \nabla_V F_H \right|^2  - \left| \tau_V^H(F_H) \right|^2 \; dV .\]
\enL

\pf First we note that from \rfE{et3} we have
\begin{align*}
\int \Gamma^V_2(f) dV &= \int \left| \nabla_H F_V\right|^2 - \mathcal{T}_1(F_V,F_H) -\mathcal{T}_2(F_H,F_H)  dV \\
&= \int \left| \nabla_H F_V\right|^2 - 2\left| \tau_V^H(F_H) \right|^2 - \aip{F_H}{ \tr[V](\nabla \tor)(F_H)}{} dV \\
&= \int \left| \nabla_H F_V\right|^2 - 2\left| \tau_V^H(F_H) \right|^2 -2 \sum\limits_k \aip{ \nabla_{U_k} F_H}{\tor(U_k,F_H)}{} dV
\end{align*}
where $U_k$ is an orthonormal frame for $VM$. Now for horizontal $E$ and vertical $U$, from \rfE{closed} we see
\begin{align*}
\aip{\nabla_E F}{U}{}^2 &= \left( \aip{\nabla_U F}{E}{} - \aip{\tor(E,U)}{F}{} \right)^2 \\
&= \aip{\nabla_U F}{E}{}^2  + \aip{\tor(E,U)}{F}{} ^2\\
& \qquad - 2  \aip{\nabla_U F}{E}{}\aip{\tor(E,U)}{F}{}.
\end{align*}
Thus summing over respective frames, we get
\[ \left| \nabla_H F_V \right|^2 = \left| \nabla_V F_H\right|^2  + \left| \tau_V^H(F_H) \right|^2  + 2\sum\limits_k \aip{\tor( \nabla_{U_k} F_H,U_k)}{F_H}{}.\]
The result then follows from torsion symmetry.

\epf

\bgL{S2} If $(M,g)$ is almost strictly normal, then for all $\nu >0$
\[ 2 \left|  S(F,F_H)  \right| \leq \nu \| \nabla_V F_H\|^2  + \frac{1}{\nu}\tau_H^V(F_H)^2 - 2 \aip{\tr \Tor(F_H)}{F}{} . \]
\enL

\pf  We compute that 
\begin{align*}
2\aip{\nabla_E F}{\tor(E,F_H)}{} &= 2  \aip{ \nabla_E F}{U}{} \aip{\tor(E,F_H)}{U}{}  \\
&= 2 \left( \aip{ \nabla_U F}{E}{} + \aip{\tor(U,E)}{F}{} \right) \aip{\tor(E,F_H)}{U}{}\\
\end{align*}
Thus summing over frames we get
\[
  2 \left|  S(F,F_H)  \right| \leq  \nu \| \nabla_V F_H\|^2  + \frac{1}{\nu}\tau_H^V(F_H)^2 -2 \aip{\tr \Tor(F_H)}{F_H}{}.
\]

\epf

\bgT{asn}
If $(M,g)$ is almost strictly normal and there are  constants $0\leq x< 1$ and  $\rho_1,\rho_2>0$ such that
\[\begin{split} \mathcal{R}^x(A,A)  &+ 2 \aip{\tr \Tor(A_H)}{A_H}{}  - 2 \left| \tau_V^H(A_H) \right| \, \left| \tau_H^V(A_H) \right|  \\ & \geq  \rho_1 \left| A_H \right|^2  + \rho_2\left| A_V \right|^2 \end{split}\]
then
\[ \lambda_1 \geq \dfrac{ \rho_1}{ (1-x) \frac{d}{d-1} + \omega}\]
where again $\omega = \kappa/\rho_2$.
\enT

\pf From \rfL{tr} and \rfL{S2} we get the alternative integrated curvature-dimension inequality
\begin{align*}
\int (1-x) \Gamma^H_2(f) &+ \nu \Gamma^V_2(f) dV \geq \int \mathcal{R}^x(A,A)  + 2 \aip{\tr \Tor(A_H)}{A_H}{}     \\
& \qquad \qquad + \frac{1-x}{d} \left( \lp f \right)^2 - \frac{1}{\nu} \left| \tau_H^V(F_H) \right|^2 - \nu \left| \tau_V^H(F_H) \right|^2  \; dV.
\end{align*}
Once again we apply this to $f \in C^\infty(M)$ with $\lp f = -\lambda f$ and $\int f^2 dV =1$. Set
\[ \nu = \frac{\rho_2}{\lambda +s}\]
with $s>0$. Then 
\begin{align*}
 (1-x) \left( \frac{d-1}{d} \right) \lambda^2  & + \frac{\lambda \rho_2}{\lambda+s} \int \Gamma^V(f) dV   \geq  \int  2\left| \tau_V^H(A_H) \right| \, \left| \tau_H^V(A_H) \right| + \rho_1 \lambda + \rho_2 \Gamma^V(f)  \\
 & \qquad - \frac{\rho_2}{s} \left| \tau_V^H(F_H) \right|^2 - \frac{\lambda+ s}{\rho_2} \left| \tau_H^V(F_H) \right|^2  \; dV.
 \end{align*}
 Now the elementary identity that for $a,b,x>0$
 \[  ax + \frac{b}{x} \geq 2 \sqrt{a}{b} \]
 with equality when $x = \sqrt{b/a}$ implies that
 \[  (1-x) \left( \frac{d-1}{d} \right) \lambda^2 + \frac{\lambda}{\rho_2} \left| \tau_H^V(F_H) \right|^2 \geq \rho_1 \lambda .\]
 The result then follows from \rfE{omega}.

\epf

\section{Examples}\label{S:EX}

\subsection{Strictly normal and $\nabla \tor$ trace-free}

If $(M,g)$ is a strictly normal sRC-manifold such that $\tr \nabla \tor (HM) =0$, then an optimal estimate is easy to obtain.

From \rfT{asn}, we immediately obtain that if 
\[ \src(A_H,A_H) \geq \rho_1 \|A_H\|^2, \qquad  \left| \tau^H_{H} (A_V) \right|^2 \geq 4\rho_2 \left|A_V\right|^2\]
then for all $0 \leq x <1$
\[ \lambda_1 \geq \dfrac{ (1-x)\rho_1 }{ (1-x) \frac{d}{d-1} + \frac{\omega}{1+3x} } = \dfrac{ \rho_1}{\frac{d}{d-1} + \frac{\omega}{(1+3x)(1-x) }}.\]
This is maximized for $x=\frac{1}{3}$, yielding
\bgE{sntf}
\lambda_1 \geq \dfrac{\rho_1}{ \frac{d}{d-1} + \frac{3\omega}{4 }}.
\enE

\subsection{Compact Lie groups}

Suppose that $G$ is a compact Lie group with Lie algebra $\mathfrak{g}$ splitting
\bgE{hv} \mathfrak{g} = \mathfrak{h} \oplus \mathfrak{v}\enE
such that $\mathfrak{v}$ is a Lie subalgebra and $\mathfrak{h}$ bracket-generates $\mathfrak{g}$ at step 2. Let $H$ be the bundle spanned by the left-invariant vector fields associated to $\mathfrak{h}$ and $V$ the corresponding bundle for $\mathfrak{v}$. Suppose also that $\mathfrak{g}$ admits a positive definite inner product that makes the split \rfE{hv} orthogonal and the operators 
\[\pi_{\mathfrak{h}} \circ \text{ad}_T \colon \mathfrak{h} \to \mathfrak{h}, \qquad \pi_{\mathfrak{v}} \circ\text{ad}_X \colon \mathfrak{v}\to \mathfrak{v}\]
skew-symmetric for all $X \in \mathfrak{h}$ and $T \in \mathfrak{v}$. (Since $G$ is compact, it admits bi-invariant metrics. If $H$ and $V$ are orthogonal with respect to a bi-invariant metric, this skew-symmetry follows immediately.) This metric is then extended to a Riemannian metric $g$ on $G$ using the correspondence between $\mathfrak{g}$ and left invariant vector fields.

For left invariant vector fields $X,Y$ in $H$ and $T,U$ in $V$, the covariant connection and the curvatures can be easily computed in terms of the Lie bracket coefficients. In the bi-invariant case: 
\[ \nabla_X Y = \frac{1}{2} [X,Y]_H, \quad \nabla_X T =[X,T]_V, \quad \nabla_TX =[T,X]_H, \quad \nabla_T U = \frac{1}{2}[T,U] . \]
Viewed as a metric extension of $g_{|H}$, the metric $g$ is strictly normal.  In the special case that $[\mathfrak{h},\mathfrak{h}]= \mathfrak{v}$, the tensor $\nabla \tor$ is guaranteed to vanish when restricted to horizontal vector fields.  With this assumption we can use \rfE{sntf}, which as we shall see in the example below greatly simplifies studying how estimates for $\mathcal{R}^x$ depend on $x$.

It should also be remarked that for any unit length left invariant vector field $X$ on a Lie group with left invariant metric, $\nabla_X X=0$. Thus the Laplace operator $\lp$ is the sum-square operator for any orthonormal frame for $HG$. Thus this methodology produces eigenvalue estimates for a wide variety of sub-elliptic sum-square operators on Lie groups.

\bgX{so4}
Let $G$ be the Lie group $SO(4)$ with Lie algebra identified with the space of skew-symmetric $4\times 4$ matrices. Denote by $M_{ij}$, the matrix with $1$ at the $(i,j)$ position, $-1$ at position $(j,i)$ and zeros elsewhere. Set $M_{12}^b =M_{12}+bM_{34}$. Define a sRC-manifold by declaring the left invariant extensions $X_1^b \sim M_{12}^b, X_2 \sim M_{13}, X_3 \sim M_{14},X_4 \sim M_{23},X_5 \sim M_{24}$ to be an orthonormal frame for $HG$ with  $T \sim M_{34}$ playing the same role for $VG$. When $b=0$, this is the standard metric on $SO(4)$. Varying $b$ represents a small perturbation of the horizontal space $HM$. The purpose of this is to illustrate the effect of the $\nabla \tor$ term within $\mathcal{R}^x$.

 The commutation table is then

\vspace{10pt}

\begin{tabular}{|c|cccccc|}
\hline
$[,]$ & $X_1^b$ &$X_2$&$X_3$&$X_4$&$X_5$&$T$\\
\hline
$X_1^b$ & 0 &$-X_4-bX_3$ &$ -X_5+bX_2$ &$X_2-bX_5 $&$ X_3+bX_4 $&0 \\
$X_2$ &$ X_4+bX_3$ &0 &$-T $&$-X_1^b+bT $&0 &$X_3$ \\
$X_3$ & $X_5-bX_2$ & $T$ &0 &0 &$-X_1^b+bT$ &$-X_2$ \\
$X_4$ & $-X_2+bX_5$ & $X_1^b-bT$ &0 &0 &$-T$ & $X_5$ \\
$X_5$ & $-X_3-bX_4$ &0 &$X_1^b-bT$ &$T$ &0 &$-X_4$ \\
$T$ & 0 &$-X_3$ &$X_2$ &$-X_5$ &$X_4$ &0 \\
\hline
\end{tabular}

\vspace{10pt}

From this it is clear that the metric condition holds and so $g$ is strictly normal.  Set $A=a^iX_i +t T$. It is also then straight-forward, if tedious, to compute that 
\begin{align*}
 \src(A_H,A_H) &= a_1^2 + \frac{3}{2} \sum\limits_{i=2}^5 a_i^2, \\
 \tr \nabla \tor(A_H) &= -a_1 b T,\\
 \left| \tau^H_{H} (A_V) \right|^2 &= 4(1+b^2)t^2,\\
 \left| \tau_H^V(A) \right|^2 &= (a_3-ba_4)^2+ (-a_2-ba_5)^2+(ba_2+a_5)^2+(ba_3-a_4)^2\\
 &\leq (1+b^2) \sum\limits_{i=2}^5 a_i^2.
 \end{align*}
From this we obtain that for all $\e>0$
 \begin{align*}
  \mathcal{R}^x(A,A) &\geq  (1-x)a_1^2 -(1+x)ba_1 t  +(1+3x)(1+b^2) t^2+ \frac{3(1-x)}{2} \sum\limits_{i=2}^5 a_i^2 \\
  & \geq  \left(1-x - \frac{\e}{2} (1+x) b \right) \|A_H\|^2  + \left( (1+3x)(1+b^2) - \frac{b(1+x)}{2\e} \right) \|A_V\|^2\\
  \kappa &= 1+b^2.
  \end{align*}
We can therefore choose $\e$ so as to obtain an estimate when
\bgE{be} \frac{(1-x)(1+3x)}{(1+x)^2} > \frac{1}{4} \frac{b^2}{1+b^2}. \enE
In particular the choice $x=0$ can always result in a (typically non-sharp) estimate.
The for all $ \frac{b}{2(1+b^2)(1+3x)} < \e < \frac{2(1-x)}{(1+x)b}$ we have
\[ \lambda_1 \geq \dfrac{ 1-x - \frac{\e}{2} (1+x) b }{\frac{4(1-x)}{5} + \frac{1+b^2}{ (1+3x)(1+b^2) - \frac{b(1+x)}{2\e} }}\]
Maximizing this in $x$ and $\e$ is a standard calculus exercise without an attractive solution.

However when $b=0$ we can use \rfE{sntf} to obtain
\[ \lambda_1 \geq \frac{20}{31} \]
%In general, we can note that \rfE{be} holds when $x=0$, to give
%Since this always holds for $x=0$, we immediately obtain that for all $\e >0$,
%\[ \lambda_1 \geq \dfrac{1- \frac{\e b}{2} } { \frac{4}{5} + \frac{ (1+b^2)}{(1+b^2)  -%\frac{b}{2\e}}}\]

\enX

\bgX{so42}
Again, let $M = SO(4)$ with the same notation as before. Thus time however, set $b=0$ and let $X_1,X_2,X_3$ be an orthonormal frame for $HM$ and $T_1=T, T_2=X_4, T_3=X_5$ be an orthonormal frame for $VM$. From the multiplication table, it is clear that $VM$ is integrable and $HM$ bracket generates at step 2. The extension is strictly normal and the horizontal Laplacian is then
\[ \lp = X_1^2 +X_2^2 +X_3^2 .\]
Repeating the calculations for $A =a^i X_i + t^k T_k$, we get
\begin{align*}
\src(A_H) &= 2\left|A_H\right|^2 = \left| \tau_H^V(A) \right|^2,  \\
 \tr \nabla \tor(A_H) &= 0,\\
 \left| \tau_H(A)\right|^2 &= 2 \left|A_V\right|^2.\\
\end{align*}
Thus we can use \rfE{sntf} with $\rho_1= 2$, $\rho_2 = \frac{1}{2}$, $\kappa =2$ and hence $\omega = 4$ to see
\[ \lambda_1 \geq \dfrac{ 2}{\frac{3}{2} + 3} = \frac{4}{9}. \]
\enX

\bgX{s3}
Now let $M=SO(3)$ with the same style of presentation as the previous example. With the same notation, we set $X^c_1$ to be the left-invariant extension of $M_{12}+cM_{13}$, with  $X_2 \sim M_{13}$ and $T \sim M_{23}$. The commutation table is then 

\vspace{10pt}

\begin{tabular}{|c|ccc|}
\hline
$[,]$ & $X^c_1$ &$X_2$&$T$\\
\hline
$X_1^c$ & 0 & $-T$ & $(1+c^2)X_2-cX_1$ \\
$X_2$ &$T $ &0  & $-X^c_1+cX_2 $ \\
$T$ & $cX_1 -(1+c^2)X_2$ & $X_1-cX_2$& 0\\
\hline
\end{tabular}

\vspace{10pt}

This time however, the metric does not meet the required conditions and the metric extension is not strictly normal unless $c=0$. It is however, almost strictly normal.

It is easy to check that $\nabla_{X_i} X_j =0$ for all $i,j=1,2$. For horizontal 
\[ \mathcal{L}_T \bp X_1 \\ X_2 \ep = \bp c & -1-c^2) \\  1 & -c \ep \bp X_1 \\ X_2 \ep .\]
Decomposing into symmetric and skew-symmetric pieces yields 
\[ \bp c & -(1+c^2) \\  1 & -c \ep = \bp 0 &  -1-c^2/2\\ 1+c^2/2 & 0  \ep - \bp -c & c^2/2 \\ c^2/2 & c \ep \]
and so
\begin{align*}
\nabla_T \bp X_1 \\ X_2 \ep &= \bp 0 &  -1-c^2/2\\ 1+c^2/2 & 0  \ep \bp X_1 \\ X_2 \ep\\
\tor\left(T, \bp X_1 \\ X_2 \ep \right) &= \bp -c & c^2/2 \\ c^2/2 & c \ep \bp X_1 \\ X_2 \ep.
\end{align*}
From this it is straightforward to compute that for $A=a X_1 +b X_2 +tT$
\begin{align*}
\mathcal{R}^x(A,A) +2\aip{\tr \Tor(A_H)}{A_H}{} &=(1-x) (a^2+b^2) + \frac{1+3x}{2} t^2 \\ & \quad  - (1+x) \left( 2abc + \frac{a^2c^2}{2} - \frac{b^2c^2}{2}\right).
\end{align*}
Now
\[ \left| \tau_H^V(A) \right|^2= a^2+b^2, \quad  \left| \tau_H^V(A) \right|^2 = \frac{a^2+b^2}{4} c^2 \left(c^2+4\right).\] 
Thus
\begin{align*}
\mathcal{R}^x(A,A) &+2\aip{\tr \Tor(A_H)}{A_H}{} -  2 \left| \tau_H^V(A) \right| \, \left| \tau_H^V(A) \right| \\
& =  \left( 1-x  - \frac{1+x}{2} c^2 - |c|\sqrt{c^2+4} \right)a^2   -2(1+x)c ab   \\
& \quad  + \left( 1-x  + \frac{1+x}{2} c^2 - |c|\sqrt{c^2+4}  \right)b^2 +\frac{1+3x}{2} t^2 \\
& \geq \left( 1-x  - \frac{3}{2} (1+x)|c| \sqrt{c^2+4} \right)(a^2+b^2) +\frac{1+3x}{2} t^2.
\end{align*}
Hence by \rfT{asn}, for all $0\leq x <1$
\[\lambda_1 \geq \dfrac{  1-x  - \frac{3}{2} (1+x)|c| \sqrt{c^2+4} }{ \frac{1-x}{2} + \frac{2}{1+3x} }.\]
Again, we are left with an elementary, if unpleasant, optimization problem in $x$ to find the best estimate. 

However, if $c=0$, then we can again use \rfE{sntf} to see that $\lambda_1 \geq \frac{1}{2}$.

\enX

\subsection{Strictly pseudoconvex pseudohermitian manifolds}

The methods described here are considerably more general than those traditionally used to study pseudohermitian manifolds. There are refinements possible in the pseudohermitian case that do not appear to have analogues in the general case. However it is still useful to see how our results apply in this case.

A pseudohermitian manifold $M^{2n+1}$ is a odd dimensional manifold equipped with a non-vanishing $1$-form $\eta$ and an endomorphism $J \colon \Ker{\eta} \to \Ker{\eta}$ with $J^2=-1$. The manifold is strictly pseudoconvex if the bilinear form $d\eta(X,JY)$ is positive definite on $\Ker{\eta}$ and hence is a sub-Riemannian metric for $\Ker{\eta}$. There is then a unique characteristic vector field $T$ such that $\eta(T)=0$, $d\eta(T,\cdot)=0$.  Setting
\[ HM = \Ker{\eta}, \qquad VM = \langle T \rangle \]
makes $M$ an sRC-manifold and we can choose a metric extension by  defining $JT=0$ and setting
\[ g(A,B) = d\eta(A,JB) + \eta(A)\eta(B).\]
It was shown in \cite{Hladky4} that the canonical connection of \rfT{Connection} is exactly the well-known Tanaka-Webster connection.  It is then easy to see that with this metric extension $M$ is totally rigid and $V$-normal. Furthermore 
\begin{align*}
\tr (\nabla \tor )(E_i) &= \sum\limits_k  \nabla_{E_k} \tor(E_i,E_K) - \tor(\nabla_{E_k} X,E_k) - \tor(X,\nabla_{E_k} E_k) \\
&=\nabla_{JE_i} T  - \aip{ \nabla_{E_k} E_i}{JE_k}{} T - \aip{E_i}{ J \nabla_{E_k} E_k}{}T \\
&= 0
\end{align*}
and
\begin{align*}
\tr \Tor(E_i)  &= \sum\limits_k \tor(E_k,\tor(E_k,E_i)) = -\tor(JE_i,T) = J \tor(E_i,T ).
\end{align*}
Hence
\begin{align*} \mathcal{R}^x(A,A)  &= (1-x) \Ric(A_H,A_H) + \frac{(1+3x) n}{2} \|A_V\|^2   \\ & \qquad -(1-x) \aip{A_H}{J \tor(A_H,T )}{}. \end{align*}
In the non-Sasakian case, the dependence of $\rho_1$ and  $\rho_2$ on $x$ is hard to determine. However, we can immediately note that 
 \[ \tau_H^V(A_H) =  \|A_H\|, \qquad \tau_V^H(A_H) = \| \tor(T,A_H) \| .\] 
 Thus if we assume that
 \bgE{spcassum}
 \begin{split}
  \src (A_H,A_H) &\geq \rho \|A_H \|^2 \\   \| \tor(T,A_H) \|\, \|A_H\| &- \aip {A_H}{J \tor(A_H,T )}{}  \leq \frac{C\rho }{2} \|A_H\|^2 
  \end{split}
  \enE
  then for all $0 \leq x <1$,
  \begin{align*} \lambda_1 &\geq \dfrac{ (1-x-C)\rho}{(1-x) \frac{2n-1}{2n} + \frac{ 2}{(1+3x) n}}  \\
  &= 2n \rho  \dfrac{ -3 x^2 +(2 -3C) x + 1 -C  }{-3(2n-1)x^2+2(2n-1)x +2n+3}. \end{align*}
  Differentiating with respect to $x$ using the quotient rule and ignoring the denominator, yields  \bgE{yuck} -2n\rho \left( 9C(2n-1) x^2  +\left( 24 +6(2n-1)C \right) x - 8  +(2n+11)C \right) . \enE
  Thus if $ C \geq \frac{8}{2n+11} $  then any critical points must have $x \leq 0$. Thus if 
  \[ \frac{8}{2n+11} \leq  C < 1\]
   we have an optimal estimate of
  \[ \lambda_1 \geq \dfrac{2n\rho}{2n+3} (1-C).\]
  However if $\rho > \frac{2n+11}{8} C$, this estimate will still hold, but the optimal estimate will occur for $x$ at the positive critical point for \rfE{yuck}, which will by necessity will be smaller than $1$.
  
 If $M$ is Sasakian, then  $C=0$. The quadratic \rfE{yuck}  collapses to $-2n\rho( 24 x - 8)$ and so the estimate is optimized at $x=1/3$. This yields
  \[ \lambda_1 \geq \frac{n \rho}{n+1}\]
  recovering \rfE{spc} at least in the torsion-free case. This last result also follows from \rfE{sntf}.

  \subsection{Twisted spheres.} This is an example of an sRC-manifold that does not fit into the framework described in \cite{BaudoinGarofalo}.  Let $M = \sn{3} \times \sn{2}$. We can view $\sn{3} \subset \rn{4}$  and describe it using coordinates $(x^1,x^2,x^3,x^4)$ and view $\sn{2} \subset \rn{3}$ using coordinates $(t^1,t^2,t^3)$ on $\rn{3}$.  We set $VM$ to be the kernel of $\pi_*$ where $\pi \colon M \to \sn{3}$ is the projection map and define vertical vector fields by
 \[ T_1 = t^2 \pd{}{t^3} - t^3\pd{}{t^2}, \quad T_2 = t^3 \pd{}{t^1} -t^1 \pd{}{t^3}, \quad T_3 = t^1 \pd{}{t^2} - t^2 \pd{}{t^1} .\] 
 We complete the sRC-structure for $M$ be declaring the following vector fields to be a global orthonormal frame for $HM$:
\begin{align*}
 X_1 &= x^2 \pd{}{x^1} - x^1 \pd{}{x^2} + x^3 \pd{}{x^4} - x^4 \pd{}{x^3} -T_1, \\
 X_2 &=  -x^4 \pd{}{x^1} + x^3 \pd{}{x^2} -x^2 \pd{}{x^3} +x^1 \pd{}{x^4} -T_2 ,\\
 X_3  &= -x^3 \pd{}{x^1} +x^1 \pd{}{x^3} - x^4\pd{}{x^2} + x^2 \pd{}{x^4} -T_3.
\end{align*}
The bracket structures are then given by
\[ [X_1,X_2]= 2X_3 -T_3, \qquad [T_1,T_2]= -T_3 \]
and
\[ [T_1,X_1]=0, \quad [T_1,X_2] = T_3, \quad [T_1,X_3]=-T_2\]
with all others determined by cyclic symmetry. A natural metric extension is defined by using the standard spherical metric on $VM$. It's then easy to compute
 \begin{align*}
 \nabla_{X_1} X_2 &= X_3, \quad \nabla_{X_2} X_1 = -X_3, \quad \tor(X_1,X_2) =T_3 ,\\
\nabla_{X_i} X_i &=0, \quad \nabla_{T_j} X_i = 0, \quad \tor(T_j,X_i) =0,  \\
\nabla_{X_j} T_i & = [X_j,T_i] 
 \end{align*} 
 with all other like terms again determined by cyclic symmetry. Thus the metric extension is strictly normal and an easy computation shows
 \begin{align*}
 R(X_1,X_2)X_2 = X_1, \qquad R(X_1,X_2)X_3 = 0 
 \end{align*}
 with cyclic symmetry and elementary properties of curvature determining the other terms. In particular, all horizontal sectional curvatures are equal to $+1$ and $\src = 2g$ on $HM$. Furthermore
\begin{align*}
 \tr \nabla \tor(X_1) & =  \nabla_{X_2} \tor(X_1,X_2) - \tor(\nabla_{X_2} X_1,X_2)  \\
 & \qquad +   \nabla_{X_3} \tor(X_1,X_3) - \tor(\nabla_{X_3} X_1,X_3)  \\
&= -\nabla_{X_2}  T_3 + \tor(X_3,X_2) + \nabla_{X_3} T_2 - \tor(X_2,X_3) \\
&= -T_1 +T_1 -T_1 +T_1 \\
&=0
 \end{align*}
 with all others similar. Thus $\tr \nabla \tor(HM) =0$. It is now an elementary exercise to show that $\kappa=2$ and
  \begin{align*}
 \mathcal{R}^x(A,A) &= 2(1-x) \left| A_H \right|^2 + \frac{1+3x}{4} \sum\limits_{i=1}^3 \aip{A}{T_i}{}^2 \\
 & \geq  2(1-x) \left| A_H \right|^2 + \frac{1+3x}{4} \left| A_V \right|^2.
 \end{align*}
 Thus for $0\leq x <1$,
 \[ \lambda_1 \geq \dfrac{ 2(1-x) }{ (1-x) \frac{2}{3} + \frac{8}{1+3x}}= \dfrac{6}{2 + \frac{24}{(1+3x)(1-x)}}.\]
 Again, this is maximized for $x=1/3$, yielding
 \[ \lambda_1 \geq \frac{3}{10}. \]

  %\begin{align*}
%\mathcal{R}^x(A,A)  \geq (1-x) \rho \|A_H\|^2 + \frac{(1+3x) n}{2} \|A_V\|^2  \end{align*}
%\[ \| A_H\| \, \|\tor(T,A_H)\| -2\aip{A_H}{J \tor(A_H,T )}{} \leq  C \left|A_H\right|^2 \]
 
 %\[ \rho_1 = (1-x) \rho-C \]
 %\[ \lambda_1 \geq  \dfrac{ (1-x) ( (1-x)-C) }{(1-x) \frac{2n-1}{2n} +  \frac{2}{(1+3x) n}} = \dfrac{ (1-x)\rho-C }{\frac{2n-1}{2n}  + \frac{ 2}{(1-x)(1+3x) n}} \]
   \bibliographystyle{plain}
  
\bibliography{References}

\end{document}